\documentclass{amsart}
\usepackage{amsfonts,hyperref}

\def\cal{\mathcal}
\def\E{\mathbf{E}} 
\def\N{\mathbb{N}} 
\def\P{\mathbf{P}} 
\def\R{\mathbb{R}} 

\theoremstyle{remark}
\newcounter{remnum}
\newtheorem{remark}{Remark}[remnum]

\newcounter{propnum}
\newtheorem{proposition}{Proposition}[propnum]

\makeatletter
\@namedef{subjclassname@2020}{\textup{2020} Mathematics Subject Classification}
\makeatother

\newcounter{todocounter}

\makeatletter
\def\@MRExtract#1 #2!{#1}
\renewcommand{\MR}[1]{
  \xdef\@MRSTRIP{\@MRExtract#1 !}
  \href{http://www.ams.org/mathscinet-getitem?mr=\@MRSTRIP}{MR\@MRSTRIP}}
\makeatother
\makeatletter
 \def\@textbottom{\vskip \z@ \@plus 5pt}
 \let\@texttop\relax
\makeatother

\begin{document}
\title[Harmonic Bernoulli trials]{On Bernoulli trials with unequal harmonic success probabilities}

\date{\today}

\author{Thierry Huillet}
\address{Thierry Huillet\newline
Laboratoire de Physique Th\'{e}orique et Mod\'{e}lisation\newline
CY Cergy Paris University, CNRS UMR-8089\newline
Site de Saint Martin, 2 avenue Adolphe-Chauvin\newline
95302 Cergy-Pontoise, France}
\email{thierry.huillet@cyu.fr}

\author{Martin M\"ohle}
\address{Martin M\"ohle\newline Fachbereich Mathematik\newline
Eberhard Karls Universit\"at T\"ubingen\newline
Auf der Morgenstelle 10\newline
72076 T\"ubingen, Germany}
\email{martin.moehle@uni-tuebingen.de}



\begin{abstract}
A Bernoulli scheme with unequal harmonic success probabilities is investigated, together with some of its natural extensions. The study includes the number of successes over some time window, the times to (between) successive successes and the time to the first success. Large sample asymptotics, statistical parameter estimation, and relations to Sibuya distributions and Yule--Simon distributions are discussed. This toy model is relevant in several applications including reliability, species sampling problems, record values breaking and random walks with disasters.
\end{abstract}

\maketitle

\section{Introduction}
Bernoulli sequences have a long history in probability and statistics with several applications, for example in reliability (Hoshino \cite{Hoshino2001}), species sampling problems (Ewens \cite{Ewens1972}, Johnson, Kotz and Balakrishnan \cite{JohnsonKotzBalakrishnan1997}, Pitman \cite{Pitman1995}), record values breaking (Neuts \cite{Neuts1967}, R\'enyi \cite{Renyi1962a,Renyi1962b}) and random walks with disasters (\cite{GoncalvesHuillet2020,Huillet2011}).

Let $w_1>0$ and $w_2\ge 0$ be two parameters. Main parts of this article deal with the two-parameter Bernoulli model, where at step $m\in\N:=\{1,2,\ldots\}$ a success occurs with probability
\begin{equation} \label{twoparamodel}
   \frac{w_1}{w+m-1},
\end{equation}
where $w:=w_1+w_2$. The particular model with $w_2=0$ or, equivalently, $w_1=w$, is extensively studied in the literature. This one-parameter model is closely related to the Ewens sampling formula \cite{Ewens1972} with mutation parameter $\theta:=w$, the Ewens fragmentation process (see, for example, Gnedin and Pitman \cite{GnedinPitman2007}) and to cycles of permutations and Poisson spacings (see, for example, Najnudel and Pitman \cite{NajnudelPitman2020}). The two-parameter model (\ref{twoparamodel}) can at least be traced back to Pitman \cite{Pitman1995} and has been studied by several authors in different context (Pitman and Yor \cite{PitmanYor1997}), with different parameter notation (see, for example, Kozubowski and Podg\'orski \cite{KozubowskiPodgorski2018}) and with different goals (see, for example, Sibuya \cite{Sibuya2014}). It is impossible to cite here all relevant literature. Research in this area is permanently ongoing; see for example the recent preprint \cite{SilvaJamshidpeyTavare2022}. For works including statistical applications, computational studies, case studies and real data applications we refer the reader to Chen and Liu  \cite{ChenLiu1997}, Hoshino \cite{Hoshino2001} and Sibuya \cite{Sibuya2014}.

Particular distributions play an important role in the context of the one- and two-parameter Bernoulli model (\ref{twoparamodel}), among them Sibuya distributions (Sibuya \cite{Sibuya1979}, Kozubowski and Podg\'orski \cite{KozubowskiPodgorski2018}) and Yule--Simon distributions (Yule \cite{Yule1925}, Simon \cite{Simon1955,Simon1960}).

One aim of this article is to review known results and to add some new insights on the two-parameter model (\ref{twoparamodel}). We furthermore extend parts of these results to a three-parameter model having an additional third parameter $\alpha\in[0,1]$. In this three-parameter Markov model, a success at step $m$ occurs with probability
\begin{equation} \label{threeparamodel}
   \frac{w_1+k\alpha}{w+m-1},
\end{equation}
where $k$ denotes the number of successes, which already occurred before step $m$. To the best of the authors knowledge, \cite{Moehle2021} is the only work where this three-parameter model implicitly occurs. Although some particular cases are mentioned in the literature (see, for example, Holst \cite{Holst2007,Holst2009} and Sibuya \cite{Sibuya2014}), the model encompassing all three parameters, $w_1$, $w_2$ and $\alpha$, appears here for the first time.

The article is organized as follows. In Section \ref{model} the two-parameter model (\ref{twoparamodel}) is introduced in detail, some basic properties of the model are discussed, possible examples and potential applications are briefly mentioned and relations to Sibuya distributions and Yule--Simon distributions are reviewed in some more detail. Section \ref{numberofsuccesses} starts with probabilistic properties of the number $S_n$ of successes in a sample of size $n$ and ends with results on parameter estimation based on data sequences of observed failures and successes. Section \ref{timestosuccessivesuccesses} contains results on the times to successive successes and Section \ref{timetofirstsuccess} on the particular time to the first success, including results on parameter estimation and hypothesis testing based on an observed sequence of such times. The remaining Sections \ref{extension}, \ref{randomwalk} and \ref{threeparameter} deal with three different extensions or modifications of the two-parameter model. Section \ref{extension} discusses a three-parameter extension where successes are more frequent. Section \ref{randomwalk} considers a related random walk with disasters. Section \ref{threeparameter} studies the more general three-parameter Markov model (\ref{threeparamodel}) for the number of successes having strong relations to (extensions of) the Chinese restaurant process and P\'olya urns.

\section{Model} \label{model}

Introduce two weights $w_1>0$ and $w_2\ge 0$, put $w:=w_1+w_2>0$, and let $I_1,I_2,\ldots$ be independent Bernoulli random variables with `harmonic' success probabilities
\begin{equation}
\P(I_m=1):=\frac{w_1}{w+m-1},\qquad m\in\N:=\{1,2,\ldots\},\label{I}
\end{equation}
decreasing inversely proportional to the number $m$ of the trial. We note the following property of Bernoulli trials with such success probabilities: the first success time $K_1^+:=\inf\{m\in\N:I_m=1\}$ is either a small number or a very large one due to power-law tails of this random variable, see (\ref{T+}) below. In words, if the $I_{m}$'s fail to take the value $1$ in the first steps, this tendency will be enhanced in the forthcoming steps resulting, for such models, in large (heavy-tailed) values of $K_1^+$. So $K_1^+$ either will take small values close to $1$ (the mode of $K_1^+$ is at $1$ with probability mass $w_1/w$ decreasing with $w_2/w_1$ if $w_2>0$) or very large values (responsible of its heavy-tailedness with tail index $w_1$): small values of $w_2/w_1$ favor early first success time while small values of $w_1$ favor late first success. So, the larger the number of steps for which no success was observed, the smaller the probability to see a success in the next step even though this probability is relatively large (harmonic decay in our case). This may be seen from the following argument:

Let $J_m:=1-I_m$, $m\in\N$, and let $M_n=\prod_{m=1}^nJ_m$. The event $M_n=1$ is realized when no success was observed till time $n$. $M_n$ is a multiplicative random walk
\[
M_{n+1}=M_nJ_{n+1},\quad M_0=1,
\]
for which the probability of a success at step $n+1$ given no success till $n$ is $\P(M_{n+1}=0\,|\,M_n=1)=\P(I_{n+1}=1)$. If $w_2=0$, then the probability $\P(M_n=1)=\prod_{m=1}^n\P(I_m=0)$ of no success by time $n\in\N$ is obviously equal to $0$, since $I_1=1$ in this case. If $w_2>0$ then this probability is equal to
\[
\prod_{m=1}^n\P(I_m=0)=\prod_{m=1}^n\frac{w_2+m-1}{w+m-1}
=\frac{\Gamma(w)\Gamma(w_2+n)}{\Gamma(w_2)\Gamma(w+n)}
\sim\frac{\Gamma(w)}{\Gamma(w_2)}\frac{1}{n^{w_1}},\quad n\to\infty,
\]
since $\Gamma(c+n)\sim n^c\Gamma(n)$ as $n\to\infty$ for any $c>0$. Thus, the probability of no success by time $n$ is small for large $n$, since $w_1>0$.

Examples of such enhancement mechanisms are
\begin{itemize}
   \item $I_m=1$ if some paper is cited the day $m$ after its publication (oversight).
   \item $I_m=1$ if some new species is discovered the day $m$ after a systematic daily sampling campaign (rareness).
   \item $I_m=1$ if some new word is used (or created) as the $m$-th word of some ongoing book (scarcity).
   \item $I_m=1$ if some individual renews its support to some political party the day (month) $m$ after its creation (weariness).
   \item Time unit increases by $1$ when some athlete attempts to improve some record previously established. $I_m=1$ if he/she succeeds at $m$-th trial: higher records become more and more difficult to break.
\end{itemize}
In several situations a success is actually a failure. Examples are
\begin{itemize}
   \item $I_m=1$ if some device breaks down the day $m$ after it was put into service (resilience).
   \item $I_m=1$ if some population collapses the day $m$ after it came to birth (resilience).
   \item $I_m=1$ if some patient contracts some illness the day $m$ after birth date (immunity).
   \item $I_m=1$ if some driver has an accident the day $m$ after obtaining his driving licence (experience).
\end{itemize}
The number $n$ of observations can be finite (possibly large though, depending on the time scale) or infinite. For instance, a typical driver only has finitely many driving days in his life (possibly randomly finite), but the attempts to break a record are potentially infinitely many.

For `harmonic' Bernoulli sequences of the form (\ref{I}) we study the number $S_n:=\sum_{m=1}^nI_m$ of successes among the first $n\in\N_0$ trials, the time $K_l^+:=\inf\{m\in\N:S_m=l\}$ of the $l$-th success, $l\in\N_0$, and the times $L_l^+:=K_l^+-K_{l-1}^+$ elapsed between successive successes, $l\in\N$, and analyse the associated Markov chains. It turns out that Sibuya distributions play an important role in this context. The two-parameter $(w_1,w_2)$-Sibuya distribution arises as the distribution of the waiting time till the first success. The shifted $(w_1,w_2)$-Sibuya distribution has many appealing properties, among them discrete self-decomposability and heavy-tailedness, \cite{KozubowskiPodgorski2018}. It includes the `bare' Sibuya distribution ($w_1+w_2=1$, see \cite{Sibuya1979}) and the Yule--Simon distribution ($w_2=1$, see \cite{Yule1925}). The case $w_2=0$ is degenerate as far as the waiting time for the first success is concerned, but it appears to make sense from the point of view of the number of successes in the Ewens species sampling problem \cite{Ewens1972}. The case $(w_1,w_2)=(1,0)$ also appears in the study of the number of record values stemming from an arbitrary independent and identically distributed (iid) sequence of observations, see \cite{Neuts1967,Renyi1962a,Renyi1962b}.

\section{Number of successes} \label{numberofsuccesses}

In this section we are mainly interested in the number $S_n:=\sum_{m=1}^nI_m$ of successes among the first $n\in\N_0$ trials. Note that $0\le S_n\le n$ for $n\in\N_0$. In particular, $S_0=0$. For (computational) results on the distribution (function) of $S_n$ and more general Poisson binomial distributions we refer the reader to Hong \cite{Hong2013} and the references therein.

In the following, $s(n,k)$, $n,k\in\N_0:=\{0,1,\ldots\}$, denote the Stirling numbers of the first kind. Recall that the unsigned Stirling numbers of the first kind $|s(n,k)|$ are characterized via $[z]_n=\sum_{k\ge 0}|s_{n,k}|z^k$, $z\in\R$, $n\in\N_0$, where $[z]_0:=1$ and $[z]_n:=z(z+1)\cdots(z+n-1)$, $n\in\N$. These numbers satisfy
the recursion $|s_{n+1,k}|=n|s_{n,k}|+|s_{n,k-1}|$ with $|s_{n,k}|=0$ for $k>n$, $|s_{n,n}|=1$ and $|s_{n,0}|=\delta_{n,0}$ (Kronecker symbol).

\subsection{The Markov chain \texorpdfstring{$(S_n,n\in\N_0)$}{}}

\indent

Clearly, $(S_n,n\in\N_0)$ is a time-inhomogeneous Markov chain with state-space $\N_0$ and transition probabilities
\begin{equation}
\P(S_{n+1}=k+1\,|\,S_n=k)=1-\P(S_{n+1}=k\,|\,S_n=k)=\frac{w_1}{w+n},
\quad n,k\in\N_0.\label{trans}
\end{equation}
Note that the probability (\ref{trans}) that the chain moves from state $k$ at time $n$ to state $k+1$ at time $n+1$ does not depend on the current state $k$. The increments
$S_n-S_{n-1}=I_n$, $n\in\N$, are independent but not identically distributed. The chain $(S_n,n\in\N_0)$ also coincides with the chain studied in the restaurant process with a cocktail bar \cite[Section 6.1]{Moehle2021} with parameters $(\alpha,\theta_1,\theta_2):=(0,w_1,w)$, where $S_n$ counts the number of occupied tables after $n$ customers have entered the restaurant. The probability generating function (pgf) $z\mapsto f_n(z):=\E(z^{S_n})$ of $S_n$ is given by
\[
f_n(z)
=\prod_{m=0}^{n-1}\frac{w_1z+w_2+m}{w+m}
=\frac{[w_1z+w_2]_n}{[w]_n}
=\frac{[w_1(z-1)+w]_n}{[w]_n},\quad z\in\R,
\]
Clearly, $f_n$ is a polynomial of degree $n$ of the form
\[
f_n(z)=\frac{1}{[w]_n}\sum_{l=0}^n|s_{n,l}|(w_1z+w_2)^l
=\frac{1}{[w]_n}\sum_{k=0}^nz^kw_1^k\sum_{l=k}^n\binom{l}{k}|s_{n,l}|w_2^{l-k}.
\]
Denoting by $[z^k]f_n(z)$ the coefficient in front of $z^k$ of $f_n$ yields
\begin{equation}
\P(S_n=k)
=[z^k]f_n(z)
=\frac{w_1^k}{[w]_n}\sum_{l=k}^n\binom{l}{k}|s_{n,l}|w_2^{l-k},
\qquad k\in\{0,\ldots,n\}.\label{pmf}
\end{equation}
With $(n)_0:=1$ and $(n)_l:=n(n-1)\cdots(n-l+1)$ for $l\in\N$, $S_n$ has the $l$-th descending factorial moment
\begin{equation}
\E((S_n)_l)
=l![(z-1)^l]f_n(z)
=\frac{w_1^l}{[w]_n}\sum_{k=l}^n(k)_l|s_{n,k}|w^{k-l},
\qquad l\in\N_0.\label{DFM}
\end{equation}
Note that $\E((S_n)_l)=0$ for $l>n$. The distribution $\pi_n(k):=\P(S_n=k)$ of $S_n$ can be recursively computed via $\pi_0(k)=\delta_{k,0}$ and
\begin{equation}
\pi_{n+1}(k)=\frac{w_1}{w+n}\pi_n(k-1)+\frac{w_2+n}{w+n}\pi_n(k),
\quad n,k\in\N_0.\label{R1}
\end{equation}
Note that $\pi_n(k)=0$ for $k\notin\{0,\ldots,n\}$. Comparing (\ref{R1}) with the recursion \cite[Theorem 1]{HsuShiue1998} $s_r(n+1,k)=s_r(n,k-1)+(n+r)s_r(n,k)$ for the generalized Stirling numbers $s_r(n,k):=S(n,k;-1,0,r)$, $n,k\in\N_0$, $r\in\R$, in the notation of \cite{HsuShiue1998} having vertical generating functions \cite[Theorem 2]{HsuShiue1998} $k!\sum_{n\ge 0}s_r(n,k)t^n/n!=(1-t)^{-r}(-\log(1-t))^k$, $r\in\R$, $k\in\N_0$, $|t|<1$, it follows that (\ref{pmf}) can be alternatively expressed in terms of these generalized Stirling numbers as
\begin{equation}
\pi_n(k)=\frac{w_1^k}{[w]_n}s_{w_2}(n,k),\qquad k\in\{0,\ldots,n\},\label{pmf2}
\end{equation}
in agreement with \cite[Eq.~(14)]{Moehle2021} for $(\alpha,\theta_1,\theta_2):=(0,w_1,w)$. Similarly, (\ref{DFM}) can be written as
\begin{equation}
   \E((S_n)_l)=\frac{w_1^l}{[w]_n}l!s_w(n,l),\qquad l\in\N_0.
\end{equation}
Introducing the superdiagonal stochastic transition matrices
\[
\Pi_n:=\left(
\begin{array}{llll}
\frac{w_2+n}{w+n} & \frac{w_1}{w+n} & 0 & \cdots\\
0 & \frac{w_2+n}{w+n} & \frac{w_1}{w+n} & 0\\
0 & 0 & \frac{w_2+n}{w+n} & \cdots\\
\vdots & \vdots & 0 & \ddots
\end{array}
\right),\qquad n\in\N_0,
\]
the distributions $\mathbf{\pi}_n:=(\pi_n(k),k\in\N_0)$ of $S_n$, $n\in\N_0$, satisfy the recursion $\mathbf{\pi}_{n+1}=\mathbf{\pi}_n\Pi_n$, $n\in\N_0$. Thus, $\mathbf{\pi}_n=\mathbf{\pi}_0\prod_{m=0}^{n-1}\Pi_m$, $n\in\N_0$, with $\mathbf{\pi}_0=(1,0,0,\ldots)$.
\begin{remark}
The pgf $f_n$ of $S_n$ has only real zeros $-(w_2+m)/w_1$, $m\in\{0,\ldots,n-1\}$. By \cite[Proposition 1]{Pitman1997}, $(\pi_n(0),\ldots,\pi_n(n))$ is a P\'{o}lya frequency sequence. Thus, the infinite matrix $M:=(\pi_n(k-l))_{k,l\in\N_0}$ (where $\pi_n(k)=0$ for $k\notin\{0,\ldots,n\}$) is totally positive of any arbitrary order, i.e., all principal minors of any arbitrary order of $M$ have nonnegative determinant.
\end{remark}
\begin{remark} (Random number of observations)
It can be natural to assume that the number of observations is finite but random (and independent of $I_1,I_2,\ldots$). In this case one has to replace $n$ by a random variable $N$ taking values in $\N$, and $\E(z^{S_N})=\sum_{n\ge 1}\E(z^{S_n})\P(N=n)$ yields the law of the number of successes over the time window $N$ with supposedly (or not) known mean $\E(N)$. For example, $N$ could be geometrically distributed $\P(N=n)=p(1-p)^{n-1}$, $n\in\N$, with parameter $p\in(0,1)$. For instance, it can be a good modeling issue to infer that there are only finitely many days in a species sampling campaign, geometrically distributed (without any further information but its mean number). The random variable $S_N$ then counts the total number of sampled species over the observation window $N$.
\end{remark}

\subsection{Special cases}

\indent

- $w=1$: $f_n(z)=\E(z^{S_n})=\frac{[w_1(z-1)+1]_n}{n!}=\frac{1}{n!}\sum_{k=0}^n|s_{n+1,k+1}|w_1^k(z-1)^k$ showing that $S_n$ has $k$-th descending factorial moment
\[
\E((S_n)_k)=k![(z-1)^k]f_n(z)=\frac{k!}{n!}|s_{n+1,k+1}|w_1^k,
\qquad k\in\N_0.
\]
Note that in that case, necessarily $w_1\in(0,1)$.

- $w_2=1$: $f_n(z)=\E(z^{S_n})=\frac{[w_1z+1]_n}{[w]_n}
=\frac{[w_1z]_{n+1}}{z[w_1]_{n+1}}
=\frac{1}{[w_1]_{n+1}}
\sum_{k=0}^{n+1}|s_{n+1,k}|w_1^kz^{k-1}$ showing that ($|s_{n,0}|=\delta_{n,0}$)
\[
\pi_n(k)=\P(S_n=k)
=\frac{|s_{n+1,k+1}|w_1^{k+1}}{[w_1]_{n+1}},\qquad
k\in\{0,\ldots,n\}.
\]

- $w_2=0$: $f_n(z)=\E(z^{S_n})
=\frac{[w_1z]_n}{[w]_n}=\frac{1}{[w_1]_{n+1}
}\sum_{k=0}^n|s_{n,k}|w_1^kz^k$ showing that
\[
\pi_n(k)=\P(S_n=k)
=\frac{|s_{n,k}|w_1^k}{[w_1]_n},\qquad k\in\{0,\ldots,n\}.
\]
If in addition $w_1=1$, then $S_n$ is the number of record values of an arbitrary iid sequence of observations appearing before $n$; \cite{Neuts1967,Renyi1962a,Renyi1962b}. In this case the law $\pi_n(k)=|s_{n,k}|/n!$, $k\in\{0,\ldots,n\}$, of $S_n$ coincides with the distribution of the number of cycles of a permutation of size $n$ chosen uniformly at random.

\subsection{Poisson approximation}

\indent

Clearly, $\mu_n:=\E(S_n)=\sum_{m=1}^n\P(I_m=1)=w_1\sum_{m=0}^{n-1}1/(w+m)=w_1\log n+O(1)$ and
\[
\sigma_n^2:={\rm Var}(S_n)=\sum_{m=1}^n\P(I_m=0)\P(I_m=1)
=w_1\sum_{m=0}^{n-1}\frac{w_2+m}{(w+m)^2}\sim w_1\log n
\]
as $n\to\infty$. The law of $S_n$ is in total variation distance close to the law of $N_n\overset{d}{\sim}\text{Poi}(\mu_n)$, (see \cite
{Pitman1997} and \cite{Yamato2017}), because $\mu_n-\sigma_n^2=w_1^2\sum_{m=0}^{n-1} 1/(w+m)^2\ll \mu_n$ (see \cite[Theorems 1 and 2]{BarbourHall1984}) with LeCam Poisson approximation of the total variation distance $d_{TV}(S_n,N_n):=\frac{1}{2}\sum_{k\ge 0}|\pi_n(k)-\mu_n^ke^{-\mu_n}/k!|$ given by (see \cite{Sevastyanov1972})
\begin{equation}
\frac{1}{32}\min (1,\mu_n^{-1}) (\mu_n-\sigma_n^2)\le d_{TV}(S_n,N_n)\le (1-e^{-\mu_n})
\frac{\mu_n-\sigma_n^2}{\mu_n}.\label{Poiapp}
\end{equation}
Therefore (and also by the Lindeberg--Feller central limit theorem), $(S_n-\mu_n)/\sigma_n\to\mathcal{N}(0,1)$ in distribution as $n\to\infty$, consistently with the fact that $(N_n-\mu_n)/\sigma_n\to\mathcal{N}(0,1)$ in distribution as $n\to\infty$. Since $\sum_{n\ge 2}\P(I_n=0)\P(I_n=1)/(\log n)^2<\infty$, it follows from well-known law of large numbers results for sums of independent, but not identically distributed random variables, that $(\log n)^{-1}\sum_{m=1}^n(I_m-\E(I_m))\to 0$ almost surely or, equivalently, that $S_n/\mu_n\to 1$ almost surely as $n\to\infty$.

\subsection{Method of moments}

\indent

Let $k:=\sum_{m=1}^ni_m$ denote the number of observed successes of a given observed sample $i=(i_1,\ldots,i_n)\in\{0,1\}^n$ of $I:=(I_1,\ldots,I_n)$.

If $w_1>0$ is known, then by the method of moments, an estimator $\widehat{w}_2$ for $w_2\ge 0$ is any solution (provided it exists) to the equation $k=\E(S_n)$ or, equivalently,
\begin{equation} \label{mm1}
   \frac{k}{w_1}=\frac{\E(S_n)}{w_1}
   =\Psi(w_1+\widehat{w}_2+n)-\Psi(w_1+\widehat{w}_2)
   =\sum_{m=0}^{n-1} \frac{1}{w_1+\widehat{w}_2+m}.
\end{equation}
For $k=0$ there is no solution $\widehat{w}_2\ge 0$ to the equation (\ref{mm1}). For $k\in\{1,\ldots,n\}$ there is a solution $\widehat{w}_2$ to the equation (\ref{mm1}) if and only if $k/w_1\le\Psi(w_1+n)-\Psi(w_1)$ and in this case the solution is unique.

If, instead, $w_2\ge 0$ is known, then, similarly, a method of moments estimator $\widehat{w}_1$ for $w_1>0$ is any solution (provided it exists) to the equation
\begin{equation} \label{mm2}
k
=\sum_{m=0}^{n-1}\frac{\widehat{w}_1}{\widehat{w}_1+w_2+m}.
\end{equation}
For $k=0$ there is no solution $\widehat{w}_1>0$ to (\ref{mm2}). For $k\in\{1,\ldots,n-1\}$ there is a unique solution $\widehat{w}_1>0$ to (\ref{mm2}), since the map $(0,\infty)\ni x\mapsto\sum_{m=0}^{n-1}x/(x+w_2+m)$ is monotone increasing with image $(0,n)$, which contains the value $k$. For $k=n$ there is no solution to (\ref{mm2}), except for the degenerate situation when $w_2=0$ and $n=1$, in which case any $\widehat{w}_1>0$ is a solution to (\ref{mm2}).

Assume now that both parameters, $w_1>0$ and $w_2\ge 0$, are unknown. Although $S_n:=\sum_{m=1}^nI_m$ admits finite moments of all orders, the following result essentially states that the standard method of moment estimator for $(w_1,w_2)$ does not exist.

\begin{proposition} \label{notexist}
   The standard method of moments estimator for $(w_1,w_2)\in(0,\infty)\times[0,\infty)$ (based on the first and second moment of $S_n$) does not exist, except for the trivial case of only $n=1$ single observed success $i_1=1$, in which case any pair $(w_1,0)$ with $w_1>0$ is a moment estimator for $(w_1,w_2)$.
\end{proposition}
\begin{proof}
   Let $i=(i_1,\ldots,i_n)\in\{0,1\}^n$ be an observed sequence of $I=(I_1,\ldots,I_n)$ having $k:=\sum_{m=1}^ni_m$ successes and let $p>0$. Since each data point $i_m$ is either $0$ or $1$, the \mbox{$p$-th} sample moment $\widehat{m}_p:=n^{-1}\sum_{m=1}^ni_m^p=k/n$ does not depend on $p$. By the standard method of moments principle, one should (try to) estimate the parameter $(w_1,w_2)$ from the two equations $\E(U_n^p)=\widehat{m}_p=k/n$, $p\in\{1,2\}$, where $U_n:=S_n/n$. These two equations imply that $\E(U_n(1-U_n))=\E(U_n)-\E(U_n^2)=\widehat{m}_1-\widehat{m}_2=k/n-k/n=0$. Thus, $U_n\in\{0,1\}$ almost surely or, equivalently, $S_n\in\{0,n\}$ almost surely. For $n\ge 2$ this is equivalent to $\P(I_1=\cdots=I_n)=1$, an obvious contradiction, since $\P(I_1\ne I_2)=w_1(2w_2+1)/(w(w+1))>0$. Thus, for $n\ge 2$ the standard method of moments estimator for $(w_1,w_2)$ does not exist. It remains to consider the case $n=1$. Then, $U_n=I_1$. Thus, the equation $k/n=\E(U_n^p)$ reduces (independently all $p>0$) to $k=\P(I_1=1)=w_1/w$. For $k=0$ there is no solution $(w_1,w_2)\in(0,\infty)\times[0,\infty)$ to this equation. For $k=1$ any pair $(w_1,0)$ with $w_1>0$ is a solution.
\end{proof}
\begin{remark}
   Let $n\in\N\setminus\{1\}$. The previous proof shows that for any sequence $I=(I_1,\ldots,I_n)$ of Bernoulli random variables, whose joint distribution depends on two parameters $w_1$ and $w_2$, the standard method of moments estimator for $(w_1,w_2)$ does not exist whenever $\P(I_1=\cdots=I_n)<1$.
\end{remark}
\begin{remark}
The result in Proposition \ref{notexist} is rather destructive. A way to circumvent this issue is the following. If the ratio $r:=\E(I_1)=w_1/w$ ($>0$) is known, then a method of moments estimator $\widehat{w}$ for $w>0$ is the solution (provided it exists) to the equation
\begin{equation} \label{mm3}
   \frac{k}{r}=\sum_{m=0}^{n-1}\frac{\widehat{w}}{\widehat{w}+m}.
\end{equation}
For $k=0$ there is no solution $\widehat{w}>0$ to (\ref{mm3}). For $k\in\{1,\ldots,n\}$ there is a solution $\widehat{w}$ to (\ref{mm3}) if and only if $k/r<n$ and in this case the solution is unique.

Note that, if the ratio $r$ is a priori not known, one may estimate $r$ in advance via a separate simple method of moments step via the number of successes divided by $n$ in an iid sample of size $n$ taken from $I_1$. This gives first an estimate for $r$, then an estimate for $w$ via (\ref{mm3}), so estimates for both, $w_1$ and $w_2$. Note that, by the central limit theorem, the estimator for $r$ is asymptotically normal with a Berry--Ess\'een rate of convergence of order $O(n^{-1/2})$.
\end{remark}

\subsection{Maximum likelihood estimation} \label{MLE_based_on_I}

\indent

As before, let $i:=(i_1,\ldots,i_n)\in\{0,1\}^n$ be an observed sample of $I=(I_1,\ldots,I_n)$ and let $k:=\sum_{m=1}^ni_m$ denote the number of observed successes. The likelihood function $L:(0,\infty)\times[0,\infty)\to[0,1]$ of the sample $i$ is given by
\begin{eqnarray} \label{joint}
   L(w_1,w_2)
   & := & \P(I=i)
   \ = \ \prod_{m=1}^n\frac{w_1^{i_m}(w_2+m-1)^{1-i_m}}{w+m-1}\nonumber\\
   & = & \frac{w_1^k\prod_{m=1}^n(w_2+m-1)^{1-i_m}}{[w]_n}.
\end{eqnarray}
This probability is not symmetric in $i_1,\ldots,i_n$ since the random variables $I_1,\ldots,I_n$ are (independent but) not identically distributed and, hence, not exchangeable. Let ${\cal L}:=\log L$ denote the log-likelihood function. Using $\partial_w\log[w]_n=\sum_{m=0}^{n-1}1/(w+m)=\Psi(w+n)-\Psi(w)$, where $\Psi:=(\log\Gamma)'=\Gamma'/\Gamma$ denotes the digamma function, the two equations $\partial_{w_j}{\cal L}(w_1,w_2)=0$, $j\in\{1,2\}$, yield \begin{equation}
\frac{k}{\widehat{w}_1}=\sum_{m=0}^{n-1}\frac{1}{\widehat{w}
+m}=\Psi(\widehat{w}+n)-\Psi(\widehat{w})
\label{MLE1}
\end{equation}
and
\begin{equation}
\sum_{m=0}^{n-1}\frac{1-i_{m+1}}{\widehat{w}_2+m}=\Psi(\widehat{w}+n)-\Psi( \widehat{w}),
\label{MLE2}
\end{equation}
where $\widehat{w}:=\widehat{w}_1+\widehat{w}_2$. Any solution $(\widehat{w}_1,\widehat{w}_2)$ to (\ref{MLE1}) and (\ref{MLE2}), provided that such a solution exists, is a maximum likelihood estimator (MLE) of $(w_1,w_2)$ based on the observed sequence $i=(i_1,\ldots,i_n)$. The Hessian matrix $H=H(\widehat{w}_1,\widehat{w}_2)$ of the log-likelihood function ${\cal L}$ at an interior point $(\widehat{w}_1,\widehat{w}_2)\in(0,\infty)^2$ of the parameter space is given by
\begin{equation} \label{hessian}
H=H(\widehat{w}_1,\widehat{w}_2)=
\left(\begin{array}{cc}
   c-a & c\\
   c & c-b
\end{array}\right)
\end{equation}
with
\[
a:=\frac{k}{\widehat{w}_1^2},\quad b:=\sum_{m=0}^{n-1}\frac{1-i_{m+1}}{(\widehat{w}_2+m)^2}
\quad\mbox{and}\quad
c:=\sum_{m=0}^{n-1}\frac{1}{(\widehat{w}+m)^2}=\Psi'(\widehat{w})-\Psi'(\widehat{w}+n).
\]
Note that $a,b\ge 0$ and $c>0$. The matrix $H$ has the characteristic polynomial $x\mapsto x^2+(a+b-2c)x+ab-ac-bc$ and, hence, the two real eigenvalues (roots of the characteristic polynomial)
\[
\lambda_{\pm}:=c-\frac{a+b}{2}\pm \sqrt{c^2+\Big(\frac{a-b}{2}\Big)^2}.
\]
Moreover, $\lambda_+<0$ (and, hence, also $\lambda_-<0$) if and only if $(a+b)c<ab$, showing that $H$ is negative definite if and only if $(a+b)c<ab$. In this case the likelihood function $L$ has a local maximum at the point $(\widehat{w}_1,\widehat{w}_2)$ in the interior $(0,\infty)^2$ of the parameter space. Otherwise there might exist a global maximum at a point $(\widehat{w}_1,\widehat{w}_2)$ with $\widehat{w}_2=0$ (Ewens model, see Section \ref{ewens}) or there exists no maximum in the parameter space $(0,\infty)\times[0,\infty)$. Recall that a MLE belonging to the interior $(0,\infty)^2$ of the parameter space exists if $(a+b)c<ab$. This condition is however not satisfied in general. For several observations $i=(i_1,\ldots,i_n)$, in particular for small sample sizes $n$, a MLE in $(0,\infty)\times[0,\infty)$ does not exist. Concrete cases illustrating these small sample size effects are provided in Table \ref{table1}.

Fortunately, if the true parameter $(w_1,w_2)$ belongs to the interior $(0,\infty)^2$ of the parameter space, then the probability that the MLE $(\widehat{w}_1(I),\widehat{w}_2(I))$ for $(w_1,w_2)$ exists and belongs to the interior $(0,\infty)^2$ of the parameter space tends to $1$ as $n\to\infty$. Moreover (see, for example, Hayashi \cite[Section 7.3, Proposition 7.8]{Hayashi2000}, Hoadley \cite{Hoadley1971} or Philippou and Roussas \cite[Section 3]{PhilippouRoussas1975}), as $n\to\infty$, asymptotic normality of $(\widehat{w}_1(I),\widehat{w}_2(I))$ holds at rate $n^{-1/2}$, that is, for all $(w_1,w_2)\in(0,\infty)^2$, $\sqrt{n}\big((\widehat{w}_1(I),\widehat{w}_2(I))-(w_1,w_2)\big)\to \mathcal{N}(0,\Sigma)$ in distribution as $n\to\infty$, where the covariance matrix $\Sigma$ of the bivariate limiting normal law $\mathcal{N}(0,\Sigma)$ is either the inverse of the expected Fisher information matrix or the inverse $J^{-1}$ of the observed information matrix $J$ evaluated at $(w_1,w_2)$. Note that $J=-H(w_1,w_2)$, where $H(w_1,w_2)$ denotes the Hessian at $(w_1,w_2)$ given by (\ref{hessian}).

If the model has two unknown parameters $(w_1,w_2)$ that have to be estimated, the first equation (\ref{MLE1}) gives $\widehat{w}_1=k/(\Psi(\widehat{w}+n)-\Psi(\widehat{w}))$ as a function of $\widehat{w}$ (and $k$) and so $\widehat{w}_2=\widehat{w}-\widehat{w}_1$ as a function of $\widehat{w}$. Plugging this expression of $\widehat{w}_2$ into the second equation (\ref{MLE2}) yields an equation in the single variable $\widehat{w}$ that can be solved from the data $i=(i_1,\ldots,i_n)$. An expression of both $\widehat{w}_1$ and $\widehat{w}_2$ then follows, provided that such a solution $(\widehat{w}_1,\widehat{w}_2)$ exists.

When $w=w_1+w_2=1$, only one parameter, say $w_1$, needs to be estimated and (\ref{MLE1}) and (\ref{MLE2}) yield
\[
\frac{k}{\widehat{w}_1}=\sum_{m=1}^n\frac{1-i_m}{m-\widehat{w}_1}
\text{ (entailing }\widehat{w}_1\in(0,1)\text{)}.
\]
When $w_2=1$ or $0$, only this first equation (\ref{MLE1}) is needed and the searched $\widehat{w}_1$, provided it exists, solves
\[
\frac{k}{\widehat{w}_1}=\sum_{m=0}^{n-1}\frac{1}{\widehat{w}_1+m+1}
\quad\text{or}\quad
\frac{k}{\widehat{w}_1}=\sum_{m=0}^{n-1}\frac{1}{\widehat{w}_1+m}.
\]
Note that $\Psi(\widehat{w}+n)-\Psi(\widehat{w})\sim\log n$ as $n\to\infty$. Thus, by (\ref{MLE1}), $\widehat{w}_1\sim k/\log n$ as $n\to\infty$ and, by (\ref{MLE2}), a large $n$ approximation for $\widehat{w}$ is the solution $w$ of the equation $\sum_{m=0}^{n-1}(1-i_{m+1})/(w+m)=\log n$.

In order to analyze the behavior of the likelihood function at the boundary of its domain it is useful to transform the parameter space one-to-one into a bounded domain, for instance by introducing the new parametrization
\begin{equation} \label{theta}
   \theta:=(\theta_1,\theta_2):=\bigg(\frac{1}{w+1},\frac{w_1}{w+1}\bigg)\in \Theta,
\end{equation}
where $\Theta:=\{\theta=(\theta_1,\theta_2):\theta_1,\theta_2>0,\theta_1+\theta_2\le 1\}$.
The Ewens model ($w_2=0$) corresponds to $\theta_1+\theta_2=1$. The models satisfying $w=w_1+w_2=1$ are equivalent to those satisfying $\theta_1=\frac{1}{2}$. The original parameters $w_1$ and $w_2$ are recovered
from $\theta$ via
\begin{equation} \label{thetainverse}
(w_1,w_2)=\bigg(\frac{\theta_2}{\theta_1},\frac{1-\theta_1-\theta_2}{\theta_1}\bigg).
\end{equation}
The maximum likelihood function in the new parametrization (\ref{theta}) continuously extends to the closure $\overline{\Theta}=\{(\theta_1,\theta_2):\theta_1,\theta_2\ge 0,\theta_1+\theta_2\le 1\}$ of $\Theta$. On this compact space $\overline{\Theta}$ a maximum, and, hence, a MLE always exists. We call such an MLE \emph{valid}, if it belongs to $\Theta$. Otherwise it is called \emph{invalid}.

Table \ref{table1} provides detailed information on the MLEs (in the closure $\overline\Theta$) for all data $i=(i_1,\ldots,i_n)\in\{0,1\}^n$ with sample sizes $n\in\{1,2,3,4\}$. For all cases listed in Table \ref{table1}, the MLE $\widehat{\theta}$ does not belong to the interior of $\Theta$. For example, for the data $i=(i_1,i_2)=(0,1)$, the likelihood function $L(w_1,w_2)=w_1w_2/(w(w+1))=\theta_2(1-\theta_1-\theta_2)/(1-\theta_1)$ takes its maximum value $\frac{1}{4}$ at the point $\widehat\theta=(\widehat\theta_1,\widehat\theta_2)=(0,\frac{1}{2})$ belonging to the boundary of $\Theta$. Note that all valid MLEs in Table \ref{table1} belong to the Ewens case ($\widehat{\theta}_1+\widehat{\theta}_2=1$).

\section{Times to successive successes} \label{timestosuccessivesuccesses}

For $l\in\N$ let $K_l^+:=\inf\{n\in\N:S_n=l\}$
be the time elapsed till the $l$-th success. Furthermore, put $K_0^+:=0$. The process $(K_l^+,l\in\N_0)$ is called the first-passage time process of the random walk $(S_n,n\in\N_0)$. Such processes have been studied extensively in the literature. We refer the reader to \cite{DenisovSakhanenkoWachtel2018} and the references therein. We have $\P(K_l^+>n)=\P(S_n<l)$ as the laws of $(K_l^+,S_n)$ are mutual inverse in the sense of inverse sampling (\cite[p.~192--194]{JohnsonKotz1977}. It follows from this, (\ref{Poiapp}), and the works \cite{Renyi1962a,Renyi1962b} (see also \cite[Proposition 1]{Moehle2021}), that
\begin{equation}
\frac{w_1\log K_l^+}{l}\overset{\text{a.s.}}{\to}1\text{ as }
l\to\infty \text{ and }\frac{w_1\log K_l^+-l}{\sqrt{l}}
\overset{d}{\to}\mathcal{N}(0,1) \text{ as }
l\to\infty,\label{LL}
\end{equation}
and the law of iterated logarithm for the $\log K_l^+$'s. And similarly for the time elapsed between contiguous successes, while replacing $K_l^+$ by $L_l^+:=K_l^+-K_{l-1}^+$ in (\ref{LL}) with the notable exception that the first almost sure convergence is now a convergence in probability \cite{Neuts1967}.

\subsection{The laws of the times to successive successes and times elapsed
between contiguous successes}

\indent

The law of $K_l^+$ is easily obtained as follows. Clearly, $\{K_l^+=n\}=\{S_{n-1}=l-1,I_n=1\}$. The independence of $S_{n-1}$ and $I_n$ thus yields
\begin{equation}
   \P(K_l^+=n)=\P(I_n=1)\P(S_{n-1}=l-1)=\frac{w_1}{w+n-1}\pi_{n-1}(l-1).
   \label{RC}
\end{equation}
Using (\ref{pmf}) the law of $K_l^+$ is therefore given by
\begin{equation}
\P(K_l^+=n)= \frac{w_1^l}{[w]_n}\sum_{k=l-1}^{n-1}\binom{k}{l-1}|s_{n-1,k}|w_2^{k-l+1},
\qquad n\ge l.
\label{PMF1}
\end{equation}
We also conclude that $L_{l+1}^+:=K_{l+1}^+-K_l^+=i$ is realized if and only if, for some $n\ge l$: $S_{n-1}=l-1$ and $I_n$ is a success and $S_{n+i-1}=l$ and $I_{n+i}$ is a success. Hence, with $i\ge 1$,
\begin{equation}
\P(L_l^+=i)=\sum_{n\ge l}
\frac{w_1}{w+n-1}\frac{w_1}{w+n+i-1}\pi_{n-1}(l-1)\pi_{n+i-1}(l),
\label{PMFL}
\end{equation}
where $\pi_n(l)$ is given by (\ref{pmf}). When $(w_1,w_2)=(1,0)$, it follows from (3) in \cite{Neuts1967}, developing problem $32$ on p.~268 in \cite{Karlin1966}, that
\[
\P(L_l^+>i)=\sum_{k=0}^i(-1)^k\binom{i}{k}(1+k)^{-l}.
\]
The law of $K_l^+$ can be obtained on a computer by launching a three-term recursion. Indeed, from (\ref{RC}), the recursion (\ref{R1}) on $\pi_n(l)$ yields a recursion for $\P(K_l^+=n)$ with $\P(K_l^+=n)=0$ if $n<l$. With $n\ge l$, this is
\begin{equation}
\P(K_{l+1}^+=n+1)=\frac{w_1}{w+n}\P(K_l^+=n) +\frac{w_2+n-1}{w+n}\P(K_{l+1}^+=n).\label{RR}
\end{equation}
Introducing the lower-triangular matrix $P=(P_{n,l})$, where $P_{n,l}:=\P(K_l^+=n)$, $l\le n$, we see that $P_{n+1,l+1}$ can be obtained from its north-west and north neighbors. With the knowledge of the first column of $P$ and its diagonal, this recursion becomes effective, starting from $P_{3,2}$ obtained from $P_{2,2}$ and $P_{2,1}$. For $n=l$, Eq.~(\ref{RR}) reduces to $\P(K_{l+1}^+=l+1)=(w_1/(w+l))\P(K_l^+=l)$, which yields the diagonal terms $\P(K_l^+=l)=\prod_{m=0}^{l-1}w_1/(w+m)$. The entries
$\P(K_1^+=n)$ of the first column of $P$ are given in (\ref{PMF+}) below.

\subsection{Markov structure of \texorpdfstring{$(K_l^+,l\in\N)$}{}}

\indent

The homogeneous Markov structure of the sequence $(K_l^+,l\in\N)$ follows from
\[
\P(K_{l+1}^+-m>n\,|\,K_l^+=m)=\prod_{k=0}^{n-1}\frac{w_2+m+k}{w+m+k}
=\frac{[w_2+m]_n}{[w+m]_n}=\prod_{k=m}^{m+n-1}\frac{w_2+k}{w+k},
\]
where $m\ge l$ and $n>0$. The random variable $L_{l+1}^+:=K_{l+1}^+-K_l^+\ge 1$ is the `time-lag' elapsed between the $l$-th and the $(l+1)$-th success. Its law depends on $K_l^+$. It is thus expected that, for each $l\ge m$, the larger $m$ is, the larger is $\P(K_{l+1}^+-m>n\,|\,K_l^+=m)$, because
\[
\frac{\P(K_{l+1}^+-(m+1)>n\,|\,K_l^+=m+1)}{\P(K_{l+1}^+-m>n\,|\,K_l^+=m)}
=\frac{w_2+m+n}{w+m+n}\frac{w+m}{w_2+m}>1.
\]
The chain $(K_l^+, l\in\N)$ therefore obeys a sort of reinforcement property. For general information on random processes with reinforcement we refer the reader to \cite{Pemantle2007}.

From Stirling's formula, $\Gamma(z+b)/\Gamma(z+a)\sim z^{b-a}$ as $z\to\infty$. For fixed $m\ll n$, for each $m\ge l$, we indeed get
\[
\P(K_{l+1}^+-m>n\,|\,K_l^+=m)
=\frac{[w_2+m]_n}{[w+m]_n}
=\frac{\Gamma(w+m)}{\Gamma(w_2+m)}(n^{-w_1}+O(n^{-(w_1+1)})),
\]
translating that, given $K_l^+=m$, the tails of $L_{l+1}^+$ have a tail index $w_1$. Given the $l$-th record occurred at $m\ll n$, the waiting time till the $(l+1)$-th has power-law tails with exponent $w_1$. Note however that the probability that $K_{l+1}^+-m=1$ is $w_1/(w+m)$ which is small only if $m\gg 1$. Introducing $c_m:=\Gamma(w+m)/\Gamma(w_2+m)$, for each $l\le m$, $c_{m+1}/c_m=(w+m)/(w_2+m)>1$ translating that the tails of $L_{l+1}^+$ get heavier as $m$ increases, but without affecting the tail index itself, only the prefactor.

With $m'>m\ge l\ge 1$, we similarly get
\begin{eqnarray*}
\P(K_{l+1}^+=m'\,|\,K_l^+=m) &=&\frac{w_1}{w+m'-1}\prod_{n=m}^{m'-2}\frac{
w_2+n}{w+n}\text{,}\\
\P(K_{l+1}^+=m') &=&\sum_{m\ge l}\P
(K_{l+1}^+=m'\,|\,K_l^+=m)\P(K_l^+=m) ,
\end{eqnarray*}
with initial condition $\P(K_1^+=m')$ given below in (\ref{PMF+}) if $w_2\ne 0$. The homogeneous Markov structure of $(K_l^+,l\in\N)$ appears more clearly, recalling $\P(K_l^+=m)$ is given by (\ref{PMF1}). For $w_2=0$ the initial condition should start with $\P(K_2^+=m')$ given in (\ref{w22}).

Introducing the excess time to $l$-th failure $K_l:=K_l^+-l\ge 0$, now a shifted random variable taking values in $\N_0$,
\begin{equation}
\P(K_{l+1}=n)
=\frac{w_1}{w+n}\P(K_l=n) +\frac{w_2+n}{w+n}\P(K_{l+1}=n-1).
\label{R2}
\end{equation}
Replacing $n$ by $n+l$ in (\ref{PMF1}) shows that the excess time $K_l$ has distribution
\[
\P(K_l=n)=\frac{w_1^l}{[w]_{n+l}}\sum_{k=l-1}^{n+l-1}\binom{k}{l-1}|s_{n+l-1,k}|w_2^{k-l+1},
\qquad l\in\N,n\in\N_0.
\]
Alternatively, the three-terms recursion (\ref{R2}) can be solved numerically using initially $\P(K_1=n)$ and observing $\P(K_l=0)=w_1^l/[w]_l$. Consequently, for all $n,n'\ge 0$,
\begin{equation}
\P(K_{l+1}>n'\,|\,K_l=n)
=\prod_{m=n+1}^{2n+n'+1}\frac{w_2+m}{w+m}  \label{Tr}
\end{equation}
and, with
\[
\P(K_{l+1}=n'\,|\,K_l=n)=\frac{w_1}{w+2n+n'+1}\prod_{m=n+1}^{2n+n'}
\frac{w_2+m}{w+m},
\]
\[
\P(K_{l+1}=k')
=\sum_{k\ge 0}\P(K_{l+1}=k'\,|\,K_l=k)\P(K_l=k),
\]
emphasizing the inhomogeneous Markov structure of $(K_l,l\in\N)$ as well. Setting $n,n'=0$ in (\ref{Tr}), we get in particular $\P(K_{l+1}=0\,|\, K_l=0)=w_1/(w+1)$.
Note that $K_l$ also represents the number of failures till the observation of the $l$-th success, a generalized version of the negative binomial distribution.

\section{Time to first success} \label{timetofirstsuccess}

If $l=1$, with $K_0^+:=0$, the distribution of the time to the first success reads ($L_1^+=K_1^+-K_0^+=K_1^+$)
\begin{equation}
\P(K_1^+>n)=[z^0]\E(z^{S_n})
=\frac{[w_2]_n}{[w]_n}
\sim\frac{\Gamma(w)}{\Gamma(w_2)}n^{-w_1},
\qquad n\to\infty.\label{T+}
\end{equation}
and
\begin{equation}
\P(K_1^+=n)=\frac{w_1}{w+n-1}\frac{[w_2]_{n-1}}{[w]_{n-1}}
=\frac{w_1}{w}\frac{[w_2]_{n-1}}{[w+1]_{n-1}},\qquad n\in\N.\label{PMF+}
\end{equation}
It is easily seen that the law of $K_1^+$ is unimodal with mode at $n=1$ having mass $w_1/w$.

Upon shifting, $K_1=K_1^+-1\ge 0$ has a generalized (heavy tailed with index $w_1$) Sibuya distribution \cite{KozubowskiPodgorski2018} with probability generating function (pgf)
\begin{equation}
\E(z^{K_1})=\frac{w_1}{w}F(1,w_2;w+1;z)\label{F}
\end{equation}
observing $(w+n)[w]_n=w[w+1]_n$, where $F:={_2}F_1$ is the Gauss hypergeometric function $F(a,b;c;z):=\sum_{n\ge 0}([a]_n[b]_n/[c]_n)(z^n/n!)$. The initial condition to the recursion (\ref{R2}) giving $\P(K_l=k)$ is
\[
\P(K_1=k)=\frac{w_1}{w+k}\frac{[w_2]_{k}}{[w]_{k}}
=\frac{w_1}{w}\frac{[w_2]_{k}}{[w+1]_{k}},\qquad k\in\N_0.
\]
\begin{remark}
(Time to first success in a $N-$Bernoulli trial with $N$ finite and Geometric
$(p)$).

In that case, $\mathcal{K}_l^+=\inf\{n\in\{1,\ldots,N\}:S_n=l\}$ and $\mathcal{K}_1^+=\inf\{m\in\{1,\ldots,N\}:I_m=1\}$. Therefore, $\mathcal{K}_1^+=\infty$ with probability $\P(S_N=0)=\E(\prod_{m=1}^{N}\P(I_m=0))$ and $\P(\mathcal{K}_1^+>n)=\P(S_n=0\,|\, N\ge n)$ with probability $\P(S_N>0)$, where
\[
\P(S_n=0\,|\,N\ge n)=\P(S_n=0)=\frac{[w_2]_n}{[w]_n}.
\]
So, if $w_2>0$, the new $\mathcal{K}_1^+$ has an atom at $\infty$ with mass
\[
\P(S_N=0)
=q\sum_{n\ge 1}\frac{[w_2]_n}{[w]_n}p^{n-1}=\frac{q}{p}[F(1,w_2;w;p)-1]
\]
translating that no success was registered before $N$.
\end{remark}

\subsection{Special cases}

\indent

- Sibuya: $w=1\Rightarrow w_1$, $w_2=1-w_1\in (0,1)$ with
$\E(z^{K_1})=w_1F(1,1-w_1;2;z)
=z^{-1}(1-(1-z)^{w_1})$, equivalently, $\P(K_1=k)=w_1[1-w_1]_{k}/(k+1)!$, $k\ge 0$.

- Yule-Simon: $w_2=1$, $w_1>0$ with $\E(z^{K_1})=
\frac{w_1}{w_1+1}F(1,1;w_1+2;z)$, equivalently, $\P(K_1=k)=w_1\frac{k!}{[w_1+1]_{k+1}}$.

- Ewens: $w_2=0$, $w_1>0$: this is a singular case for which $\P(K_1^+=n)=\delta_{n,1}$.

In view of $F(a,b;c;z)=F(b,a;c;z)$, the Yule--Simon distribution with $a=b=1$ and $c=w_1+2$ is the only one in the class (\ref{F}) to be identifiable (different parameters yield different distributions).

\subsection{Falling factorial moments of \texorpdfstring{$K_1$}{}}

\indent

$K_1^+=K_1+1$ is an important random variable if one considers that the first occurrence of a success may lead to a stop of some ongoing process.

With $a=1$, $b=w_2$, $c=w+1$, $i$ integer, using the special values and differential identities
\begin{eqnarray*}
F(a,b;c;1) &=&\frac{[c-a]_{a}}{[c-a-b]
_{a}}\\
\frac{{\rm d}^i}{{\rm d}z^i}F(a,b;c;z) &=&\frac{[a]
_{i}[b]_{i}}{[c]_{i}}F(a+i,b+i;c+i;z) ,
\end{eqnarray*}
evaluated at $z=1$, with $(K_1)_i=K_1(K_1-1)\cdots(K_1-i+1)$, when $i<w_1$, we get the descending $i$-th factorial moments of $K_1$ as
\[
\E[(K_1)_i]=\varphi^{(i)}(1)
=\frac{i![w_2]_{i}}{[w_1-i]_{i}},\qquad i<w_1,
\]
where $\varphi(z):=\E(z^{K_1})=\frac{w_1}{w}F(1,w_2;w+1;z)$. In particular, if $w_1>1$, $\E(K_1)=w_2/(w_1-1)<\infty$ and, if $w_1>2$,
\[
{\rm Var}(K_1)=\varphi''(1)+\varphi'(1) -(\varphi'(1))^2
=\frac{w_1(w-1)\E(K_1)}{(w_1-1)(w_1-2)}<\infty.
\]
Overdispersion holds. The mean $\E(K_1^+)=\frac{w-1}{w_1-1}>1$ and the variance ${\rm Var}(K_1^+)={\rm Var}(K_1)$ of $K_1^+$ (if they exist) may be used to
estimate $(w_1,w_2)$ by the method of moments provided empirical values of these quantities are available.

\subsection{Maximum likelihood estimation \texorpdfstring{of $(w_1,w_2)$ from $K_1^+$}{}}

\label{MLE}

\indent

The maximum likelihood estimator (MLE) discussed in Section \ref{MLE_based_on_I} is based on the independent but not identically distributed random variables $I_1,\ldots,I_n$. Instead, we are now going to estimate $(w_1,w_2)$ from a sample $(k_1,\ldots,k_n)\in\N^n$ taken from $n\in\N$ independent copies $K_1^+(1),\ldots,K_1^+(n)$ of the time $K_1^+$ of the first success. This section is kept short, since similar results on method of moments estimators (MME) and MLE of $(w_1,w_2)$ based on random samples taken from a generalized Sibuya distribution are provided in \cite[Section 6]{KozubowskiPodgorski2018}. However, we provide some further details pointing to the direction of hypothesis testing discussed in Section \ref{testing} in more detail. Clearly,
\[
\P(K_1^+(1)=k_1,\ldots,K_1^+(n)=k_n)
=w_1^n\prod_{\ell=1}^n\frac{[w_2]_{k_\ell-1}}{[w]_{k_\ell}}.
\]
Considering $\partial_{w_j}\log\P(K_1^+(1)=k_1,\ldots,K_1^+(n)=k_n)=0$ for $j\in\{1,2\}$ yields a MLE $(\widehat{w}_1,\widehat{w}_2)$ for $(w_1,w_2)$ based on the histogram of the sample $(k_1,\ldots,k_n)$ of observed times of first success, provided that a solution $(\widehat{w}_1,\widehat{w}_2)$ to the following equations (\ref{first}) and (\ref{second}) exists. With $\widehat{w}=\widehat{w}_1+\widehat{w}_2$, we get
\begin{equation} \label{first}
   \frac{n}{\widehat{w}_1}
   =\sum_{\ell=1}^n\big(\Psi(\widehat{w}+k_\ell)-\Psi(\widehat{w})\big)
\end{equation}
and
\begin{equation} \label{second}
   \sum_{\ell=1}^n\big(\Psi(\widehat{w}_2+k_\ell-1)-\Psi(\widehat{w}_2)\big)=\sum_{\ell=1}^n \big(\Psi(\widehat{w}+k_\ell)-\Psi(\widehat{w})\big).
\end{equation}
The Hessian matrix $H=H(\widehat{w}_1,\widehat{w}_2)$ of the log-likelihood function at an interior point $(\widehat{w}_1,\widehat{w}_2)\in(0,\infty)^2$ of the parameter space has again the form
\begin{equation} \label{hessian2}
   H=\left(\begin{array}{cc}c-a&c\\c&c-b\end{array}\right)
\end{equation}
as in (\ref{hessian}), but with different parameters
\[
a:=\frac{n}{\widehat{w}_1^2},\quad
b:=\sum_{\ell=1}^n\big(\Psi'(\widehat{w}_2)-\Psi'(\widehat{w}_2+k_\ell-1)\big)\quad\mbox{and}\quad
c:=\sum_{\ell=1}^n\big(\Psi'(\widehat{w})-\Psi'(\widehat{w}+k_\ell)\big).
\]
The matrix $H$ is negative definite if and only if $(a+b)c<ab$. In this case the likelihood function has a local maximum at $(\widehat{w}_1,\widehat{w}_2)$. Otherwise there might exist a global maximum at the border $\partial ((0,\infty]\times[0,\infty))=(\{0,\infty\}\times[0,\infty])\cup([0,\infty]\times\{0,\infty\})$ of the parameter range. In the latter case the MLE does not necessarily need to belong to the original parameter space $(0,\infty)\times[0,\infty)$. If the MLE is local, then the
first equation (\ref{first}) gives $\widehat{w}_1$ as a function of $\widehat{w}$ (and the data) and so $\widehat{w}_2=\widehat{w}-\widehat{w}_1$ as a function of $\widehat{w}$. Plugging this expression of $\widehat{w}_2$ into the second equation (\ref{second}) yields an equation in the single variable $\widehat{w}$ that can be solved from the data. A separate expression of both $\widehat{w}_1$ and $\widehat{w}_2$ then follows. In the alternative parametrization $(\beta,\theta):=(1/w,w_1/w)\in(0,\infty)\times(0,1]$, asymptotic normality of the corresponding estimator is shown in \cite[Theorem 1]{KozubowskiPodgorski2018}, together with an expression of the Fisher information matrix. More details on this asymptotic normality are provided in Section \ref{testing}.

\subsection{Hypothesis testing} \label{testing}

\indent

A comprehensive discussion of hypothesis testing is beyond the scope of this work. We exemplary consider a typical two-sided test situation. For some given $(w_1^0,w_2^0)\in(0,\infty)^2$ we would like to test the hypothesis $H_0:(w_1,w_2)=(w_1^0,w_2^0)$ against the alternative $H_1:(w_1,w_2)\ne(w_1^0,w_2^0)$. Intuitively, $H_0$ is rejected if some suitable estimator $(\widehat{w}_1,\widehat{w}_2)$ for $(w_1,w_2)$ is ``far away'' from $(w_1^0,w_2^0)$. We exemplary use here the MLE $(\widehat{w}_1,\widehat{w}_2)=(\widehat{w}_1(k),\widehat{w}_2(k))$ introduced in the previous Section \ref{MLE} based on a sample $k:=(k_1,\ldots,k_n)\in\N^n$ taken from $n\in\N$ independent copies $K_1^+(1),\ldots,K_1^+(n)$ of the time $K_1^+$ of the first success. We furthermore define $K:=(K_1^+(1),\ldots,K_1^+(n))$ for convenience. A reasonable test $\phi_n:\N^n\to\{0,1\}$ for $H_0$ against $H_1$ is thus of the form
\begin{equation} \label{test}
   \phi_n(k)=1_{\{\|(\widehat{w}_1(k),\widehat{w}_2(k))-(w_1^0,w_2^0)\|>c_n\}},\qquad k\in\N^n,
\end{equation}
where $\|.\|$ denotes the usual Euclidian norm in $\R^2$ and $c_n>0$ is a constant which has to be determined from the constrain that the probability
\begin{eqnarray*}
   \P_{(w_1^0,w_2^0)}(\phi_n=1)
   & = & \P_{(w_1^0,w_2^0)}\big(
            \|(\widehat{w}_1(k),\widehat{w}_2(k))-(w_1^0,w_2^0)\|>c_n
         \big)\\
   & = & \P\big(\|(\widehat{w}_1(K),\widehat{w}_2(K))-(w_1^0,w_2^0)\|>c_n\big)
\end{eqnarray*}
for a type $1$ error should not exceed a pre-specified significance level $\alpha\in(0,1)$. It is readily checked that the test $\phi_n$ based on the sample $k=(k_1,\ldots,k_n)$ has the $p$-value
\begin{equation} \label{pvalue}
   \widehat{\alpha}_n(k)
   =1-F_n\big(\sqrt{n}\,\|(\widehat{w}_1(k),\widehat{w}_2(k))-(w_1^0,w_2^0)\|\big),
\end{equation}
where $F_n$ denotes the distribution function of $\sqrt{n}\,\|(\widehat{w}_1(K),\widehat{w}_2(K))-(w_1^0,w_2^0)\|$.

A large sample size approximative test $\phi$ is obtained as follows. Since the likelihood function is infinitely differentiable on the interior of its domain, it follows that the standard conditions for the asymptotic efficiency of MLEs are satisfied. Thus (see, for example, Lehmann \cite[p.~429, Theorem 4.1]{Lehmann1983} or Lehmann and Casella \cite[p.~449, Theorem 3.10]{LehmannCasella1998}), under the hypothesis $H_0$, for any $(w_1^0,w_2^0)$ belonging to the interior $(0,\infty)^2$ of the parameter space, the asymptotic normality
\begin{equation} \label{asynormal}
   \sqrt{n}\big((\widehat{w}_1(K),\widehat{w}_2(K))-(w_1^0,w_2^0)\big)\to \mathcal{N}(0,\Sigma)
\end{equation}
in distribution as $n\to\infty$ holds with covariance matrix $\Sigma:=I^{-1}$, where the Fisher information matrix $I=-H$ is the negative of the Hessian $H=H(w_1^0,w_2^0)$ (see (\ref{hessian2})) evaluated at $(w_1^0,w_2^0)$. Thus, a large sample size approximation $c$ for the constant $c_n$ is the solution to the equation $\P(\|X\|>\sqrt{n}c)=\alpha$ or, equivalently, $c=q_{1-\alpha}/\sqrt{n}$, where $X=(X_1,X_2)$ is an $\R^2$-valued random variable with distribution $\mathcal{N}(0,\Sigma)$ and $q_{1-\alpha}$ denotes the $(1-\alpha)$-quantile of the distribution of $\|X\|$. Note that the distribution of $X$ depends on $\Sigma$ and, hence, on $(w_1^0,w_2^0)$. The corresponding approximative test has the $p$-value
\begin{equation} \label{approxpvalue}
   \widehat{\alpha}(k)
   =1-F\big(\sqrt{n}\,\|(\widehat{w}_1(k),\widehat{w}_2(k))-(w_1^0,w_2^0)\|\big),
\end{equation}
where $F$ denotes the distribution function of $\|X\|$. Since $\|X\|^2=X_1^2+X_2^2$ has the same distribution as $\lambda_1 Z_1^2+\lambda_2Z_2^2$, where $\lambda_1$ and $\lambda_2$ are the eigenvalues of $\Sigma$ and $Z_1$ and $Z_2$ are independent and both standard normal distributed, it follows that $\|X\|^2$ has Laplace transform $s\mapsto((1+2\lambda_1s)(1+2\lambda_2s))^{-1/2}$, $s\ge 0$, and density
\begin{eqnarray*}
   x
   & \mapsto & \frac{1}{2\pi\sqrt{\lambda_1\lambda_2}}
         \int_0^x\frac{1}{\sqrt{(x-y)y}}
         \exp\Big(-\frac{x-y}{2\lambda_1}-\frac{y}{2\lambda_2}\Big){\rm d}y\\
   & = & \frac{1}{2\sqrt{\lambda_1\lambda_2}}
         \exp\Big(-\frac{\lambda_1+\lambda_2}{4\lambda_1\lambda_2}x\Big)
         I_0\Big(\frac{\lambda_1-\lambda_2}{4\lambda_1\lambda_2}x\Big),
         \qquad x>0,
\end{eqnarray*}
where $I_0(t)=\sum_{k=0}^\infty t^{2k}/(2^{2k}(k!)^2)$ denotes the modified Bessel function of the first kind (of order $0$) satisfying the differential equation $x^2y''+xy'-x^2y=0$. If $\lambda_1=\lambda_2$ then $\|X\|^2$ is exponentially distributed with mean $2\lambda_1$. In general, $\|X\|^2$ has mean $\lambda_1+\lambda_2$ and variance $2(\lambda_1^2+\lambda_2^2)$. For general information on (representations for the densities and distribution functions of) quadratic forms of multivariate normal random variables we refer the reader to the book of Mathai and Provost \cite{MathaiProvost1992}.

The error of the asymptotic test $\phi$ compared to the exact test $\phi_n$ can be controlled as follows. For $n\in\N$ the quantity $r_n:=\sup_{x\in\R}|F_n(x)-F(x)|$ can be viewed as a measure of the rate or speed of the convergence. The asymptotic normality (\ref{asynormal}) together with the continuity of the norm $\|.\|$ and of the distribution function $F$ of $\|X\|$ clearly imply that $r_n\to 0$ as $n\to\infty$. Thus, by (\ref{pvalue}) and (\ref{approxpvalue}), the distance of the $p$-values of the exact and the asymptotic test is bounded by $|\widehat\alpha_n(k)-\widehat\alpha(k)|\le r_n$ uniformly for all observed data $k=(k_1,\ldots,k_n)$. We conjecture a typical Berry--Ess\'een rate of convergence of order $r_n=O(n^{-1/2})$. A rigorous proof of this conjecture would probably require advanced techniques going beyond the scope of this work.

One-sided test situations can be handled analogously. We leave the details to the interested reader.

\subsection{The Ewens case \texorpdfstring{$w_2=0$}{}} \label{ewens}

\indent

In a sampling problem from a Poisson--Dirichlet partition $\text{PD}(\theta)$ of the unit interval modeling species abundances, the law of the number $S_n=\sum_{m=1}^nI_m$ of distinct sampled species for a size $n$ uniform sample obeys (\ref{R1}), \cite{Ewens1972}, \cite{ArratiaBarbourTavare1992} and \cite{Yamato2017}, with $w_1=\theta$, $w_2=0$ and $S_1=1$, corresponding to $K_1^+=1$. Because sampling is modeled as uniform throws on a partition of the unit interval, necessarily on day $n=1$, a new species is sampled but new species with smaller abundance become increasingly unlikely to be subsequently sampled. The $\text{PD}(\theta)$ partition of the unit interval has countably many pieces, so the sampling process potentially never stops. Here $K_l^+$ ($l\ge 2$) is the sample size till $l$ new species have been sampled with, from (\ref{PMF1})
\begin{equation}
\P(K_l^+=n)
=w_1[z^{l-1}]\frac{[w_1z]_{n-1}}{[w_1]_n}
=w_1^{l-1}\frac{|s_{n-1,l-1}|}{[w_1+1]_{n-1}},\qquad n\ge l.
\label{w2=0}
\end{equation}
This distribution seems to be new. Note the resulting `vertical' identity for the $|s_{n,l}|$'s: $\sum_{n\ge l}\frac{|s_{n,l}|}{[w_1+1]_n}=w_1^{-l}$ for all $w_1>0$.

The random variable $K_2^+$ is the time to second non-trivial discovery of a new species (after $K_1^+=1$), with, recalling $|s_{n-1,1}|=(n-2)!$,
\begin{equation}
\P(K_2^+=n)
=w_1\frac{(n-2)!}{[w_1+1]_{n-1}},\qquad n\ge 2,  \label{w22}
\end{equation}
reducing to $\P(K_2^+=n)=1/(n(n-1))$ when $w_1=1$. With $K_2^+-1\ge 0$ the time
elapsed since $K_1^+=1$, we thus have
\[
\E(z^{K_2^+-1})=w_1\int_0^z\frac{F(1,1;w_1+1;t)-1}{t}{\rm d}t.
\]
The above theory applies to this fundamental Ewens model. Given $S_n=k$, the probability to discover a new species at time $n+1$ is $w_1/(w_1+n)$, decreasing inversely proportional to $n$ and independently of $k$. Recall from (\ref{MLE1}) that the MLE $\widehat{w}_1$ for $w_1$ is characterized by $k/\widehat{w}_1=\Psi(\widehat{w}_1+n)-\Psi(\widehat{w}_1)$ and hence only depends on $k=i_1+\cdots+i_n$. See \cite[p.~41, Eq. (3.7.7)]{TavareZeitouni2004}.

\section{An extension of the harmonic Bernoulli trial} \label{extension}

With $\alpha>0$, consider the inhomogeneous Bernoulli trial with $\P(I_m=1)=w_1/(w+(m-1)^\alpha)$, $m\in\N$.

For $\alpha\in(0,1)$ the successful events are more frequent than for $\alpha=1$. Then, $\mu_n:=\E(S_n)=\sum_{m=1}^n\P(I_m=1)=w_1\sum_{m=0}^{n-1}1/(w+m^\alpha)\sim
\frac{w_1}{1-\alpha}n^{1-\alpha}$ as $n\to\infty$ and
\[
\sigma_n^2:={\rm Var}(S_n)=\sum_{m=1}^n
\P(I_m=1)\P(I_m=0)
=w_1\sum_{m=0}^{n-1}\frac{w_2+m^\alpha}{(w+m^\alpha)^2}
\sim\frac{w_1}{1-\alpha}n^{1-\alpha}
\]
and the law of $S_n$ is close in the sense of total variation distance to $P_n\overset{d}{\sim}{\rm Poi}(\mu_n)$ for this new $\mu_n$ now growing algebraically with $n$.

Clearly also,
\begin{equation}
\frac{w_1(K_l^+)^{1-\alpha}}{l(1-\alpha)}\overset{\text{a.s.}}{\to}1
\text{ as }l\to\infty \text{ and }
\frac{w_1(K_l^+)^{1-\alpha}/(1-\alpha)-l}{\sqrt{l}}\overset{d}{\to}\mathcal{N}( 0,1)
\text{ as }l\to\infty.  \label{LL2}
\end{equation}
The time to the $l$-th success occurs much sooner than when $\alpha=1$.

If $\alpha>1$, then $S_n$ converges in distribution to a Poisson random variable with finite mean $\mu_\infty:=\lim_{n\to\infty}\mu_n=w_1\sum_{m=0}^\infty 1/(w+m^\alpha)$.

\section{A related random walk with disasters} \label{randomwalk}

Bernoulli trials with unequal harmonic success probabilities are also relevant in the context of growth-collapse random walks with disasters. Discrete-time integral-valued growth-collapse processes where long periods of linear growth alternate with rare catastrophic events occur in a large variety of systems. A collapse or catastrophic event is when the size of some population shrinks by a random number of units, not exceeding the current system's size. A total disaster is when the size of the system shrinks instantaneously to zero (a massive extinction event). Disastrous growth-collapse models occur as models for population growth subject to rare catastrophic extinction events.

A one-parameter version of such discrete-time models was investigated in \cite{Huillet2011}. Here, holding probabilities were allowed (with some probability the system's size can be left unchanged) and pure reflection at the origin was assumed (once in state zero, the system's size grows by one unit with probability $1$). Whenever zero is a reflection/absorption barrier, pomp periods will alternate with periods of scarcity. We herewith focus on discrete-time disastrous growth-collapse models with no holding probability and with zero either standing for a reflection or an absorption barrier. The probabilities of either growth or disastrous events will be chosen to be dependent on the current state as in the Bernoulli model with harmonic success probabilities, and this will favor large populations in the long run.

With $\alpha>0$, define $q_n:=w_1/(w+n^\alpha)$ and $p_n:=1-q_n$, $n\in\N_0$. With $(U_m,m\in\N)$ an iid sequence of uniforms,
\begin{equation}
N_{m+1}:=(N_m+1) \mathbf{1}(U_{m+1}\le p_{N_m}),\qquad N_0\ge 0,\label{MCD}
\end{equation}
defines a time-homogeneous Markov chain $(N_m,m\in\N_0)$ that moves from state $n$ to state $n+1$ with probability $p_n$ or is sent from state $n$ to state $0$ with probability $q_n$ (a disaster event).

The transition matrix of this Markov chain with state-space $\N_0$ is
\[
P=\left(
\begin{array}{cccccc}
q_0 & p_0 & & & & \cdots\\
q_1 & 0 & p_1 & & & \cdots\\
\vdots & \vdots & \ddots & \ddots & & \cdots\\
q_n & 0 & \cdots & 0 & p_n & \cdots\\
\vdots & \vdots & & & \ddots & \ddots
\end{array}
\right).
\]
Let us distinguish two cases.

\vspace{2mm}

\textbf{Case 1.} Assume that $w_2=0$. In this case state $0$ is absorbing. Let $n\in\N$. The probability $\P(N_m\to\infty\,|\,N_0=n)=\prod_{m\ge n}p_m$ is equal to $0$ if and only if $\sum_{m\ge n}q_m=\infty$, which in turn holds if and only if $\alpha\le 1$. Thus, for $\alpha\le 1$ the chain $(N_m,m\in\N_0)$, started from state $N_0\equiv n$, will eventually go extinct. For $\alpha>1$ the chain, started from state $n$, will tend to infinity with probability $\prod_{m\ge n}p_m>0$ and go extinct with complementary probability $1-\prod_{m\ge n}p_m$. The extinction time $\tau_{n,0}:=\inf\{m\in\N_0:N_m=0,N_0=n\}$ has pgf $\E(z^{\tau_{n,0}})=\sum_{m\ge n}q_mz^m\prod_{k=n}^{m-1}p_k$,
$|z|<1$, and $\tau_{n,0}$ takes the value $\infty$ with probability $\prod_{m\ge n}p_m$ being strictly positive if and only if $\alpha>1$.

\vspace{2mm}

\textbf{Case 2.} Assume that $w_2>0$. Then state $0$ is reflecting and all states are communicating since $w_1>0$ by assumption. The chain $(N_m,m\in\N_0)$ is hence irreducible and obviously aperiodic. This is a small variation of a Markov chain whose salient statistical features were studied in \cite{GoncalvesHuillet2020}. From the study in \cite{GoncalvesHuillet2020} we conclude that:
\begin{itemize}
   \item For $\alpha>1$ the chain is transient. After a finite number of returns to $0$ (excursions) the chain drifts to infinity.
   \item For $\alpha<1$ the chain is positive recurrent with invariant probability measure $\pi_n=\pi_0\prod_{k=0}^{n-1}p_k$, $n\in\N_0$, where the normalizing constant $\pi_0$ is determined by $\sum_{n=0}^\infty\pi_n=1$.
   \item For $\alpha=1$ (critical case) the chain is null-recurrent if $0<w_1\le 1$ and positive recurrent if $w_1>1$. For the latter case $w_1>1$ the invariant probability measure is given by $\pi_n=\pi_0[w_2]_n/[w]_n$, $n\in\N_0$, with normalizing constant $\pi_0:=(w_1-1)/(w-1)$, having heavy tails with index $w_1>1$.
\end{itemize}
In the recurrent case ($\alpha\le 1$) the sample paths of $(N_m,m\in\N_0)$, started at $N_0=0$, are made of iid excursions through state $0$. The first excursion has length $L_1^+$ and height $L_1^+-1$, where $L_1^+:=\inf\{m\in\N:N_m=0,N_0=0\}$ is the time elapsed till the first disaster. Clearly, in the positive recurrent case ($\alpha<1$ or $\alpha=1$ and $w_1\le 1$) the invariant probability measure has the general form $\pi_n=\P(L_1^+>n)/\E(L_1^+)$, $n\in\N_0$.

With $(L_i^+-1,i\in\N)$ iid copies of the first excursion height $L_1^+-1$, of interest for the control of overcrowding are the random variables
\[
T_1(n) :=\inf\{m\in\N:N_m>n\,|\,N_0=n_0\}
\text{ and }\inf\{i\in\N:\max_{j\in\{1,\ldots,i\}}(L_j^+-1)>n\},
\]
corresponding to the first (overcrossing) time the chain $N_m$ exceeds $n$ given $N_0=n_0<n$ and the number of the corresponding excursion.

Let $P_{(n)}$ be the truncated upper-left corner with size $(n+1,n+1)$ of the full irreducible transition matrix $P$ of $N_m$ (its north-west part). With $\mathbf{1}'=(1,\ldots,1)$ and $\mathbf{e}_{n_0}'=(0,\ldots,0,1,0,\ldots,0)$ transpose row vectors with size $n+1$ (with $1$ in position $n_0+1$ for $\mathbf{e}_{n_0}'$), it follows from \cite[Propositions 11 and 12]{GoncalvesHuillet2021} that
\begin{equation}
   \P_{n_0}(T_1(n)>l)=\mathbf{e}_{n_0}'P_{(n)}^l\mathbf{1},\label{TAIL}
\end{equation}
where $P_{(n)}^l$ is the $l$-th power of $P_{(n)}$. $\P(T_1(n)>l)=1$ for $l\in\{1,\ldots,n-n_0\}$. At this time $T_1(n)$, the state of the chain $N_m$ is $n+1$ because the overshoot can only be $1$. So $T_1(n)$ has geometric tails with decay-rate parameter the spectral radius of $P_{(n)}$ and
\[
\E_{n_0}(T_1(n))=\mathbf{e}_{n_0}'(I-P_{(n)})^{-1}\mathbf{1}.
\]
Clearly, given $N_0=n_0<n$, with $N_l^*=\max_{m\le l}N_m$ the extremal process of $N_m$, the events $N_l^*\le n$ and $T_1(n)>l$ coincide, so (\ref{TAIL}) also gives the marginal law $\P_{n_0}(N_l^*\le n)$ of $N_l^*$.

The extremal chain $N_l^*$ only grows (by one unit) at the record times
$R_k:=\inf\{r\in\N:r>R_{k-1},N_r>N_{R_{k-1}}\}$ of $N_m$.

\section{A more general Markov model for the number of successes} \label{threeparameter}
As before, let $w_1>0$ and $w_2\ge 0$ and define $w:=w_1+w_2$. A more general model can be introduced by taking an additional third parameter $\alpha\in[0,1]$ and assuming that the number $S_n$ of successes forms a Markov chain $(S_n,n\in\N_0)$ satisfying $S_0=0$ and
\[
\P(S_{n+1}=k+1\,|\,S_n=k)=1-\P(S_{n+1}=k\,|\, S_n=k):=\frac{w_1+k\alpha}{w+n},\qquad n\in\N_0.
\]
In this case $S_n$ coincides with the number of occupied tables in the (Chinese) restaurant process with a cocktail bar \cite{Moehle2021} after $n$ customers have entered the restaurant. For $\alpha=0$ we are back to the model studied before. For $\alpha>0$ the transition probabilities of the random walk $(S_n,n\in\N_0)$ now depend not only on the time $n$ but also on the current state $S_n=k$. The distribution of $S_n$ can be expressed as (see \cite[Eq.~(14)]{Moehle2021})
\[
\P(S_n=k)=\frac{[w_1|\alpha]_k}{[w]_n}S(n,k;-1,-\alpha,w_2),
\qquad k\in\{0,\ldots,n\},
\]
where $[w_1|\alpha]_0:=1$, $[w_1|\alpha]_k:=\prod_{i=0}^{k-1}(w_1+i\alpha)$ for $k\in\N$ and $S(n,k;-1,-\alpha,w_2)$ denote the generalized Stirling numbers in the notation of Hsu and Shiue \cite{HsuShiue1998}, which can be calculated as follows.  For $\alpha=0$ it follows from (\ref{pmf}) that $S(n,k;-1,0,w_2)=\sum_{l=k}^n\binom{l}{k}w_2^{l-k}|s_{n,l}|$, $k\in\{0,\ldots,n\}$. For $\alpha\ne 0$, the Dobi\'nski-type formula \cite[Theorem 4]{HsuShiue1998} yields
\[
S(n,k;-1,-\alpha,w_2)=\frac{1}{k!\alpha^k}\sum_{l=0}^k(-1)^l
\binom{k}{l}[w_2-l\alpha]_n,\qquad k\in\{0,\ldots,n\}.
\]
Note that $\P(S_n=0)=[w_2]_n/[w]_n$ does not depend on $\alpha\in[0,1]$. In particular, for any $n\in\N$, $\P(S_n=0)=0$ if and only if $w_2=0$. Formulas for the moments of $S_n$ are provided in \cite[Section 6.1]{Moehle2021} for $\alpha=0$ and in \cite[Corollary 1]{Moehle2021} for $\alpha>0$. The behavior of $S_n$ for $\alpha>0$ differs substantially from the case $\alpha=0$. For $\alpha>0$, as $n\to\infty$, $S_n/n^\alpha$ converges almost surely and in $L^p$ for any $p>0$ to a limiting random variable being three-parameter $(\alpha,\beta,\gamma)$-Mittag--Leffler distributed, where $\beta:=w$ and $\gamma:=w_1/\alpha$, see \cite[Theorem 3]{Moehle2021}. We refer the reader to \cite[Section 7]{Moehle2021} for further details on the three-parameter Mittag--Leffler distribution ${\rm ML}(\alpha,\beta,\gamma)$. For $\alpha=1$ the limiting distribution ${\rm ML}(1,w,w_1)=\beta(w_1,w_2)$ is the beta distribution with parameters $w_1$ and $w_2$, in agreement with well-known results for standard P\'olya urns.

If $w_2=0$ then $S_n$ counts the number of distinct species in a sample of size $n$ taken from Pitman and Yor's \cite{PitmanYor1997} two-parameter stick-breaking ${\rm PD}(\alpha,w_1)$-partition of the unit interval, extending the Ewens case. We refer the reader to Chapter 3 of Pitman's lecture notes \cite{Pitman2006} for further information on the two-parameter model and to Yamato and Sibuya \cite{YamatoSibuya2000} and Yamato, Sibuya and Nomachi \cite{YamatoSibuyaNomachi2001} for some further related works.

Assume now that $w_2>0$. In this case $S_n$ may no longer be seen, stricto sensu, as the number of new species in a sample of size $n$ taken from a partition of the unit interval. However (see \cite[Theorem 2]{Moehle2021}), $S_n$ is the number of new species (excluding a `fictitious species' $0$ with beta distributed `abundance' $B_0\stackrel{d}{=}\beta(w_2,w_1)$) in a sample of size $n$ drawn from a kind of three-parameter Poisson--Dirichlet partition ${\rm PD}(\alpha,w_1,w_2):=(B_0,(1-B_0)\text{PD}(\alpha ,w_1))$, where $B_0$ is independent of ${\rm PD}(\alpha,w_1)$.

Note that $K_1^+:=\inf\{n\in\N:S_n=1\}$ has distribution
\[
\P(K_1^+=n)=\P(S_{n-1}=0)
\P(S_n=1\,|\,S_{n-1}=0)=w_1\frac{[w_2]_{n-1}}{[w]_n},\quad n\in\N,
\]
so that $S_n^+:=S_{n+K_1^+-1}$ (with $S_1^+=1$) coincides (in law) with the number of new species from a PD$(\alpha,w_1)$-partition of the unit interval. Whenever a sample hits the `fictitious species' $0$, sampling simply fails to draw any new species: this event thus represents the possibility of a failure of the sampling process from scratch. The probability that in a sample of size $n$ there are $n_0$ failure events clearly is the beta binomial probability mass function $\binom{n}{n_0}[w_2]_{n_0}[w_1]_{n-n_0}/[w]_n$, $n_0\in\{0,\ldots,n\}$. If $\alpha=0$ then $S_n$ is the number of new species (excluding the `fictitious species' $0$ with `abundance' $B_0$) in a sample of size $n$ drawn from the partition ${\rm PD}(0,w_1,w_2)=(B_0,(1-B_0){\rm PD}(0,w_1))$, extending the Ewens case.

Let $n_0\in\N_0$ and $n_1,\ldots,n_k\in\N$ and put $n:=n_0+\cdots+n_k$. Note that
\begin{eqnarray}
   &   & \hspace{-5mm}\P(S_n=k,N_n(0)=n_0,N_n(1)=n_1,\ldots,N_n(k)=n_k)\nonumber\\
   &   & \hspace{5mm}=\ n!\frac{[w_1|\alpha]_k}{[w]_n}\frac{[w_2]_{n_0}}{n_0!}
   \prod_{l=1}^k \frac{[1-\alpha]_{n_l-1}}{(n_l-1)!\sum_{j=l}^kn_j}
   \label{dtg}
\end{eqnarray}
is the joint distribution that there are $n_0$ visits to the reservoir set with size $B_0$ (accounting for early failure events of the sampling process, or missed samples) and $S_n=k$ distinct visited species in order of appearance with positive sample sizes $n_1,\ldots,n_k$ not in the reservoir. For $w_2=0$, (\ref{dtg}) reduces to the two-parameter Donnelly--Tavar\'e--Griffiths distribution ${\rm DTG}(w_1,\alpha)$ (see \cite[Theorem 1]{YamatoSibuyaNomachi2001})
\begin{equation}
   \P(S_n=k,N_n(1)=n_1,\ldots,N_n(k)=n_k)
   \ =\ n!\frac{[w_1|\alpha]_k}{[w_1]_n}
        \prod_{l=1}^k \frac{[1-\alpha]_{n_l-1}}{(n_l-1)!\sum_{j=l}^kn_j}.
        \label{dtg2}
\end{equation}
For $\alpha=0$, (\ref{dtg2}) reduces to
\begin{equation}
   \P(S_n=k,N_n(1)=n_1,\ldots,N_n(k)=n_k)
   \ =\ n!\frac{w_1^k}{[w_1]_n}\prod_{l=1}^k \frac{1}{\sum_{j=l}^kn_j},
   \label{dt}
\end{equation}
which is \cite[Eq.~(1)]{Yamato1997} with $\alpha$ there replaced by $w_1$. The Stirling numbers $s(n,k)$ of the first kind satisfy the elementary relation (see, for example, \cite[Eq.~(4.3)]{DonnellyTavare1986})
\[
\sum_{\substack{n_1,\ldots,n_k\in\N\\n_1+\cdots+n_k=n}}\frac{n!}{(n_1+\cdots+n_k)\cdots(n_{k-1}+n_k)n_k}\ =\ |s(n,k)|,\quad n\in\N, k\in\{1,\ldots,n\}.
\]
Hence, summing the right hand side of (\ref{dt}) over all $k\in\{1,\ldots,n\}$ and $n_1,\ldots,n_k\in\N$ satisfying $n_1+\cdots+n_k=n$ yields the value
\[
\frac{1}{[w_1]_n}\sum_{k=1}^n w_1^k \sum_{\substack{n_1,\ldots,n_k\in\N\\n_1+\cdots+n_k=n}}
\frac{n!}{\prod_{l=1}^k\sum_{j=l}^kn_j}
\ =\ \frac{1}{[w_1]_n} \sum_{k=1}^n w_1^k |s(n,k)|
\ =\ 1,
\]
showing that (the right hand side of) (\ref{dt}) is indeed a probability distribution.

Summing (\ref{dtg}) over all $n_1,\ldots,n_k\in\N$ with $n_1+\cdots+n_k=n-n_0$, the joint probability that, in a sample of size $n$, there are $S_n=k$ new sampled species and $n_0\le n$ visits to the `fictitious species' is
\[
\P(N_n(0)=n_0,S_n=k)
\ =\ \binom{n}{n_0}\frac{[w_2]_{n_0}[w_1|\alpha]_k}
{[w]_n}S(n-n_0,k;-1,-\alpha,0).
\]
Observing that
$\sum_{k=0}^{n-n_0}S(n-n_0,k;-1,-\alpha,0)[w_1|\alpha]_k=[w_1]_{n-n_0}$, the probability that, in a sample of size $n$, there are $n_0\le n$ visits to the `fictitious species' is thus the beta-binomial probability $\P(N_n(0)=n_0)=\binom{n}{n_0}[w_2]_{n_0}[w_1]_{n-n_0}/[w]_n$, in agreement with the explanations above.

Given $N_n(0)=n_0$, one may estimate $(w_1,\alpha)$ either using the MLE on top of p.~506 of Hoshino \cite{Hoshino2001} (with $w_1$ substituted to $\theta$ there), or the moment estimator \cite[Eqs.~(22) and (23)]{Hoshino2001} based on singleton (unique) and doubleton observables (with $n-n_0$ substituted to $n$). Hoshino \cite{Hoshino2001} confronted his estimators to the field of microdata disclosures, while Sibuya \cite{Sibuya2014} applied the MLE estimator to trawl data arising from fishing industries. The original main field of application is species abundance problems arising in statistical ecology.

For an observed value of $n_0$ given a sample of size $n$, a first moment
estimator of the ratio $w_2/w$ is $n_0/n$, completing the estimation of $w_2$. This results from the fact that $N_n(0)$ is Beta-binomial distributed. In summary, this gives a method to estimate all three parameters $w_1$, $w_2$ and $\alpha$, first $(w_1,\alpha)$ and then $w_2$.

$\bullet$ MLE: This results from the first Ewens--Pitman sampling formula ($N_n(0)$ being Beta-binomial distributed):
\begin{eqnarray*}
   &   & \hspace{-15mm}\P(N_n(0)=n_0,N_n(1)=n_1,\ldots,N_n(k)=n_k,S_n=k)\\
   & = & \frac{n!}{k!}\frac{[w_1|\alpha]_k}{[w]_n}\frac{[w_2]_{n_0}}{n_0!}
         \prod_{\ell=1}^k\frac{[1-\alpha]_{n_\ell-1}}{n_\ell!}\\
   & = & \binom{n}{n_0}\frac{[w_1]_{n-n_0}[w_2]_{n_0}}{[w]_n}
         \frac{(n-n_0)!}{k!}\frac{[w_1|\alpha]_k}{[w_1]_{n-n_0}}
         \prod_{\ell=1}^k\frac{[1-\alpha]_{n_\ell-1}}{n_\ell!}
\end{eqnarray*}
and, hence,
\begin{eqnarray*}
   &   & \hspace{-15mm}\P(N_n(1)=n_1,\ldots,N_n(k)=n_k,S_n=k\,|\,N_n(0)=n_0)\\
   & = & \frac{(n-n_0)!}{k!}\frac{[w_1|\alpha]_k}{[w_1]_{n-n_0}}
         \prod_{\ell=1}^k\frac{[1-\alpha]_{n_\ell-1}}{n_\ell!},
\end{eqnarray*}
where to the right-hand side we recognize the joint distribution of the $N_n(\ell)'$s when sampling is from ${\rm PD}(w_1,\alpha)$, observed in an arbitrary way (independently of the sampling mechanism) with sample size $\sum_{\ell=1}^kn_\ell=n-n_0$.

$\bullet$ MME: If $w_2\ge 0$, with $n_0\ge 0$, $\sum_{m=1}^{n-n_0}ma_m=n-n_0$ and $\sum_{m=1}^{n-n_0}a_m=k$, the second alternative Ewens--Pitman sampling formula, accounting for repetitions of occupancies, reads
\begin{eqnarray} \label{jlaw}
   &   & \hspace{-15mm}
         \P(N_n(0)=n_0,A_n(1)=a_1,\ldots,A_n(n-n_0)=a_{n-n_0},S_n=k)\nonumber\\
   & = & n!\frac{[w_1|\alpha]_k}{[w]_n}\frac{[w_2]_{n_0}}{n_0!}
         \prod_{m=1}^{n-n_0}\frac{[1-\alpha]_{m-1}^{a_m}}{m!^{a_m}a_m!}.
\end{eqnarray}
Here, $A_n(m)$ is the number of boxes (species) visited $m$ times by the $n$-sample and $S_n=n-N_n(0)$ is the number of visited boxes.

Joint moments of the random variables $A_n(1),\ldots,A_n(n-n_0)$ are obtained as follows. With $i:=\sum_{m=1}^{n-n_0}i_m$ and $j:=\sum_{m=1}^{n-n_0}mi_m$,
given $N_n(0)=n_0$, the joint falling factorial moments of the $A_n(m)$'s are given by (following (4.1) in \cite{YamatoSibuya2000})
\begin{eqnarray} \label{jmom}
   &   & \hspace{-15mm}
         \E\Bigg(\prod_{m=1}^{n-n_0}(A_n(m))_{i_m}\Bigg| N_n(0)=n_0\Bigg)\nonumber\\
   & = & (n-n_0)_j
         \frac{[w_1|\alpha]_i[w_1+i\alpha]_{n-n_0-j}}{[w_1]_{n-n_0}}
         \prod_{m=1}^{n-n_0}\left(\frac{[1-\alpha]_{m-1}}{m!}\right)^{i_m},
\end{eqnarray}
giving access to the marginals when all integers $i_m$s are $0$ except one
and also to the two-point correlations Corr$(A_n(m),A_n(m'))$, $m'\ne m$.

\section*{acknowledgment}
T.~Huillet acknowledges partial support from the Chair `\textit{Mod\'elisation math\'ematique et biodiversit\'e'} of Veolia-Ecole Polytechnique-MNHN-Fondation X and from the labex MME-DII Center of Excellence (\textit{Mod\`eles math\'ematiques et \'economiques de la dynamique, de l'incertitude et des interactions}, ANR-11-LABX-0023-01 project). This work was also funded by CY Initiative of Excellence (grant `\textit{Investissements d'Avenir}' ANR-16-IDEX-0008), Project `EcoDep' PSI-AAP2020-0000000013.

\vfill\eject

\section*{Tables}

\begin{table}[ht] 

\centering
   \begin{tabular}{lllllc}
   $n$ & $k$ & data $i$
      & MLE $(\widehat{\theta}_1,\widehat{\theta}_2)$ & maximum value & valid\\
   \hline
   $1$ & $0$ & $(0)$
      & $(t,0)$, $t\in[0,1]$ & $1$ & no\\
   $1$ & $1$ & $(1)$
      & $(t,1-t)$, $t\in[0,1]$ & $1$ & no\\
   \hline
   $2$ & $0$ & $(0,0)$
      & $(t,0)$, $t\in[0,1]$ & $1$ & no\\
   $2$ & $1$ & $(0,1)$
      & $(0,\frac{1}{2})$ & $\frac{1}{4}$ & no\\
   $2$ & $1$ & $(1,0)$
      & $(1,0)$ & $1$ & no\\
   $2$ & $2$ & $(1,1)$
      & $(0,1)$ & $1$ & no\\
   \hline
   $3$ & $0$ & $(0,0,0)$
      & $(t,0)$, $t\in[0,1]$ & $1$ & no\\
   $3$ & $1$ & $(0,0,1)$
      & $(0,\frac{1}{3})$ & $\frac{2}{27}\approx 0.07407$ & no\\
   $3$ & $1$ & $(0,1,0)$
      & $(0,\frac{1}{3})$ & $\frac{2}{27}\approx 0.07407$ & no\\
   $3$ & $1$ & $(1,0,0)$
      & $(1,0)$ & $1$ & no\\
   $3$ & $2$ & $(0,1,1)$
      & $(0,\frac{2}{3})$ & $\frac{4}{27}\approx 0.14815$ & no\\
   $3$ & $2$ & $(1,0,1)$
      & $(\sqrt{2}-1,2-\sqrt{2})$ & $3-2\sqrt{2}\approx 0.17157$ & yes\\
   $3$ & $2$ & $(1,1,0)$
      & $(\sqrt{2}-1,2-\sqrt{2})$ & $6-4\sqrt{2}\approx 0.34315$ & yes\\
   $3$ & $3$ & $(1,1,1)$
      & $(0,1)$ & $1$ & no\\
   \hline
   $4$ & $0$ & $(0,0,0,0)$
      & $(0,t)$, $t\in[0,1]$ & $1$ & no\\
   $4$ & $1$ & $(0,0,0,1)$
      & $(0,\frac{1}{4})$ & $\frac{27}{256}$ & no\\
   $4$ & $1$ & $(0,0,1,0)$
      & $(0,\frac{1}{4})$ & $\frac{27}{256}$ & no\\
   $4$ & $1$ & $(0,1,0,0)$
      & $(0,\frac{1}{4})$ & $\frac{27}{256}$ & no\\
   $4$ & $1$ & $(1,0,0,0)$
      & $(1,0)$ & $1$ & no\\
   $4$ & $2$ & $(0,0,1,1)$
      & $(0,\frac{1}{2})$ & $\frac{1}{16}=0.0625$ & no\\
   $4$ & $2$ & $(0,1,0,1)$
      & $(0,\frac{1}{2})$ & $\frac{1}{16}=0.0625$ & no\\
   $4$ & $2$ & $(0,1,1,0)$
      & $(0,\frac{1}{2})$ & $\frac{1}{16}=0.0625$ & no\\
   $4$ & $2$ & $(1,0,0,1)$
      & $\approx(0.53209,0.46791)$ & $\approx 0.08378$ & yes\\
   $4$ & $2$ & $(1,0,1,0)$
      & $\approx(0.53209,0.46791)$ & $\approx 0.12567$ & yes\\
   $4$ & $2$ & $(1,1,0,0)$
      & $\approx(0.53209,0.46791)$ & $\approx 0.25133$ & yes\\
   $4$ & $3$ & $(0,1,1,1)$
      & $(0,\frac{3}{4})$ & $\frac{27}{256}\approx 0.10547$ & no\\
   $4$ & $3$ & $(1,0,1,1)$
      & $(0,\frac{3}{4})$ & $\frac{27}{256}\approx 0.10547$ & no\\
   $4$ & $3$ & $(1,1,0,1)$
      & $\approx(0.20980,0.79020)$ & $\approx 0.15256$ & yes\\
   $4$ & $3$ & $(1,1,1,0)$
      & $\approx(0.20980,0.79020)$ & $\approx 0.22883$ & yes\\
   $4$ & $4$ & $(1,1,1,1)$
      & $(0,1)$ & $1$ & no\\
   \hline
   \end{tabular}

   \vspace{2mm}

   \caption{MLE $\widehat{\theta}\in\overline{\Theta}$ in the parametrization (\ref{theta}) for all data sequences $i=(i_1,\ldots,i_n)\in\{0,1\}^n$ with sample size $n\in\{1,2,3,4\}$}
   \label{table1}
\end{table}

\vspace{5mm}

\section*{Conflict of interest statement}
On behalf of all authors the corresponding author states that there is no conflict of interest.

\vfill\eject

\end{document}